\documentclass[12pt]{amsart}
\usepackage{amsmath, amsthm}
\usepackage{amsmath,amssymb,amsthm}
\usepackage{amsfonts}
\usepackage[mathscr]{eucal}
\newtheorem{theorem}{Theorem}[section]

\newtheorem{corollary}[theorem]{Corollary}

\theoremstyle{definition}

\theoremstyle{remark}

\numberwithin{equation}{section}

\newcommand{\RR}[1]{\mathbb{#1}}

\newcommand{\rd}{{\mathbb R^d}}

\begin{document}

\title{\bf Symmetric $\alpha$-stable subordinators and Cauchy problems}

\author{Erkan Nane}
\address{Erkan Nane, Department Statistics and Probability,
Michigan State University, East Lansing, MI 48823}
\email{nane@stt.msu.edu}
\urladdr{http://www.stt.msu.edu/$\sim$nane}


\begin{abstract}
 
 We survey the results in Nane (E. Nane, \emph{Higher order PDE's and iterated processes},  Trans. American Math. Soc. (to appear)) and Baeumer, Meerschaert, and Nane (B. Baeumer, M.M. Meerschaert and E. Nane, \emph{Brownian subordinators and fractional Cauchy problems:} Submitted (2007)) which deal with PDE connection of some iterated processes, and obtain a new probabilistic proof of the equivalence of the higher order PDE's   and fractional in time PDE's.
\end{abstract}


\maketitle

\section{\bf Introduction }

In recent years, starting with the articles of Burdzy
\cite{burdzy1, burdzy2}, researchers have shown interest in iterated
processes in which one changes the time parameter with
one-dimensional Brownian motion.

To define {\bf iterated  Brownian motion} $Z_{t}$, due to
Burdzy \cite{burdzy1}, started at $z \in \RR{R}$,
 let $X_{t}^{+}$, $X_{t}^{-}$ and $Y_{t}$
 be three  independent
one-dimensional Brownian motions, all started at $0$. {\bf Two-sided
Brownian motion} is defined to be
\[ X_{t}=\left\{ \begin{array}{ll}
X_{t}^{+}, &t\geq 0\\
X_{(-t)}^{-}, &t<0.
\end{array}
\right. \] Then iterated Brownian motion started at $z \in \RR{R}$
is
$ Z_{t}=z+X(Y_{t}),\ \ \    t\geq 0.$


\subsection{\bf BM versus IBM}
This process is not Markovian or Gaussian but, it has many properties analogous to those of
Brownian motion. We list a few;

 IBM $Z_{t}$ has stationary (but not independent) increments, and is a
{\bf self-similar process} of index $1/4$.
{\bf Laws of the iterated logarithm (LIL)} holds:
usual LIL by Burdzy \cite{burdzy1} shows
$$
\limsup_{t\to\infty}\frac{Z(t)}{t^{1/4}(\log \log
(1/t))^{3/4}}=\frac{2^{5/4}}{3^{3/4}} \ \ \  a.s.
$$
 Chung-type LIL is obtained by Khoshnevisan and Lewis \cite{koslew} and Hu et al.
\cite{hu}.
 Khoshnevisan
and Lewis \cite{klewis} extended results of Burdzy \cite{burdzy2}, to develop a
{\bf stochastic calculus} for iterated Brownian motion.
 Burdzy and Khosnevisan \cite{bukh} showed that IBM can be used to model diffusion in a crack.
 Local times of this process was  studied by Burdzy
and Khosnevisan \cite{burdzy-khos},  Cs\'{a}ki, Cs\"{o}rg\"{o}, F\"{o}ldes,
and R\'{e}v\'{e}sz \cite{CCFR}, Shi and Yor \cite{shi-yor}, Xiao \cite{xiao}, and Hu \cite{hu-2}.
 Ba\~{n}uelos and DeBlassie \cite{bandeb} studied the
{\bf distribution of exit place} for iterated Brownian motion
in cones.
 DeBlassie \cite{deblassie} studied the lifetime asymptotics of iterated
Brownian motion in cones and Bounded domains. Nane \cite{nane, nane2, nane5, nane6}
extended some of the results of DeBlassie. 

\subsection{\bf  PDE connection}

In addition to the above properties of IBM, there is an interesting
connection between iterated Brownian motion
 and the {\bf biharmonic operator $\Delta ^{2}$}. Allouba and Zheng \cite{allouba1} show that if we replace the
outer process $X(t)$ in the definition of iterated Brownian motion  with a continuous Markov process with  $L_x$ as the semigroup generator, then
$u(t,x)=E_{x}[f(Z_{t})]:=E[f(Z_{t})|Z_0=x]$
solves the Cauchy initial value problem
\begin{equation}\label{ibm-pde00}
\frac{\partial}{\partial t}u(t,x) =
\frac{{L_x}f(x)}{\sqrt{\pi t}}+ {L_x}^{2}u(t,x); \quad u(0,x) = f(x), \ t>0,  x\in \rd .
\end{equation}
When $Z_t$ is an iterated Brownian motion, this was also obtained by DeBlassie \cite{deblassie} by a different method.
Let $Z^1_t=X(|Y_t|)$, then Allouba and Zheng \cite{allouba2} shows $E_{x}[f(Z_{t})]=E_{x}[f(Z^1_{t})]$.

\subsection{\bf Fractional Cauchy problems}

Nigmatullin \cite{nigmatullin} gave a physical derivation of the fractional kinetic
equation for some special $\beta$
\begin{equation}\label{frac-derivative-0}
\frac{\partial^\beta}{\partial t^\beta}u(t,x) =
{L_x}u(t,x); \quad u(0,x) = f(x)
\end{equation}
where $0< \beta <1$ and $L_x$ is the generator of
some continuous Markov process $X_0(t)$ started at $x=0$. Here $\partial^{\beta} g(t)/\partial t^\beta $
is the Caputo fractional derivative in time, which can be defined as the
inverse Laplace transform of $s^{\beta}\tilde{g}(s)-s^{\beta -1}g(0)$, with
$\tilde{g}(s)=\int_{0}^{\infty}e^{-st}g(t)dt$ the usual Laplace transform.

Mathematical study of equation (\ref{frac-derivative-0}) was initiated by \cite{sch-wyss, koch1, koch2}. The existence and uniqueness of solutions to equation (\ref{frac-derivative-0}) was proved in \cite{koch1, koch2}.
This equation was also used by Zaslavsky \cite{zaslavsky} for Hamiltonian chaos.


Baeumer and Meerschaert \cite{fracCauchy} and Meerschaert and Scheffler \cite{limitCTRW}
 show that the fractional Cauchy problem (\ref{frac-derivative-0}) is related to a certain
 class of subordinated stochastic processes;
 take $D_t$ to be the stable subordinator,
 a L\'evy process with strictly increasing sample paths such that $E[e^{-sD_t}]=e^{-ts^\beta}$,
  see for example Bertoin \cite{bertoin}.  Define the inverse or hitting time or first passage time process
\begin{equation}\label{Edef}
E_t=\inf\{x>0: D(x)>t\} .
\end{equation}
The subordinated process $Z_t=X_0(E_{t})$ occurs as the scaling limit of a
 continuous time random walk (also called a renewal reward process), in which iid
 random jumps are separated by iid positive waiting times (Meerschaert and Scheffler (2004)\cite{limitCTRW}).
Theorem 3.1 in Baeumer and Meerschaert \cite{fracCauchy} shows that, in the
case $p(t,x)=T(t)f(x)$ is a bounded continuous semigroup on a
Banach space, the formula
$$
u(t,x)=\int_0^\infty p((t/s)^\beta,x) g_\beta(s)\,ds
=
\frac{t}{\beta}\int_{0}^{\infty}p(x,s)g_{\beta}(\frac{t}{s^{1/\beta}})s^{-1/\beta
-1}ds \nonumber 
$$
yields a solution to the   fractional Cauchy problem (\ref{frac-derivative-0}).
  Here, $g_\beta(t)$ is the smooth
density of the stable subordinator, with $\tilde g_\beta(s)=\int_0^\infty
e^{-st}g_\beta(t)\,dt=e^{-s^\beta}$.

\section{\bf Brownian subordiantors and  fractional Cauchy problems}
We give a probabilistic proof of the following theorem. A variation of this result was realized by Orsingher and Benghin \cite{OB} for a version of iterated Brownian motion.
\begin{theorem}[Baeumer, Meerschaert and Nane (2007)\cite{nane-m-b}]\label{equiv-pde}
Let $L_{x}$ be the generator of a Markov semigroup $T(t)f(x)=E_{x}[f(X_{t})]$,
and take $f\in D(L_{x})$ the domain of the generator.  Then, both the higher order Cauchy problem (\ref{ibm-pde00})
 and the fractional Cauchy problem (\ref{frac-derivative-0})
 with $\beta=1/2$,
have the same solution
\begin{eqnarray}
u(t,x)& = &E_x[f(Z_t)]
= \frac{2}{\sqrt{4\pi t}}\int_{0}^{\infty}T(s)f(x)\exp\left(-\frac{s^2}{4t}\right)ds. \label{ACPsoln1}
\end{eqnarray}
\end{theorem}
\begin{proof}
$E_t$ is the inverse of a $1-1/\alpha$ stable subordinator. $E_t$ then is the local time of symmetric stable process of index $\alpha$. In the case $\alpha=2$, local time of Brownian motion is the same as $\sup_{0<s<t} B_s$. On the other hand, $\sup_{0<s<t} B_s$ and $|B_t|$  are same in distribution by the reflection principle. Hence, $E_t$ and $|B_t|$ have the same one-dimensional distributions, implying the result of the theorem.
\end{proof}

We obtain the following corollary of our theorem

\begin{corollary}[Baeumer, Meerschaert and Nane (2007)\cite{nane-m-b}] 
For any continuous Markov process
$X(t)$, both the Brownian-time subordinated process $X(|Y_t|)$ and
the process $X(E_t)$ subordinated to the inverse $1/2$-stable
subordinator have the same one-dimensional distributions.  Hence
they are both stochastic solutions to the fractional Cauchy
problem (\ref{frac-derivative-0}), or equivalently, to the higher
order Cauchy problem \eqref{ibm-pde00}.
\end{corollary}
In contrast to the previous corollary, we have
\begin{theorem}
Let $Y$ be a symmetric stable process of index $1<\alpha <2$, and $E_t$ is the inverse of a stable subordinator of index $1-1/\alpha$.
The processes $X(E_t)$ and $X(|Y_t|)$ do not have same one-dimensional distribution.
\end{theorem}
\begin{proof}Let $L_1^0$ be the local time at $x=0$. $L_1^0$ has the same one-dimenional distributions as $E_t$.
Lemma 1 in Hawkes \cite{hawkes} implies that
\begin{equation}\label{local-tail}
P[L_1^0>\lambda]\sim C_1\lambda^{-\alpha/2}\exp(-C_{\alpha h}\lambda^{\alpha}).
\end{equation}
 Proposition 4 in Bertoin \cite{bertoin} shows
\begin{equation}\label{stable-tail}
 P[Y_1>u]\sim P[\sup_{0\leq s \leq 1}Y_s >u]\sim cu^{-\alpha}.
\end{equation}
 
 The results in equations (\ref{local-tail}) and (\ref{stable-tail})establish  that in the case $Y_t$ is a symmetric stable process of index $\alpha <2$, $|Y_t|$ and $E_t$ do not have the same one-dimensional distributions.
\end{proof}

When the outer proces is   L\'evy process we have uniqueness of the solutions in Theorem \ref{equiv-pde}.
 The proof  relies on a Laplace-Fourier transform argument.

\begin{theorem}[Baeumer, Meerschaert and Nane (2007)\cite{nane-m-b}]\label{unique-sol}
Suppose that $X(t)=x+X_0(t)$ where $X_0(t)$ is a L\'evy process starting at zero.  If $L_x$ is the generator of the semigroup $T(t)f(x)=E_x[(f(X_t))]$ on $L^1(\rd)$, then for any $f\in D(L_{x})$, both the initial value problem \eqref{ibm-pde00}, and the fractional Cauchy problem \eqref{frac-derivative-0} with $\beta=1/2$, have the same unique solution given by \eqref{ACPsoln1}.
\end{theorem}

An easy extension of the argument for Theorem \ref{unique-sol} shows that, under the same conditions, for any $n=2,3,4,\ldots$ both the Cauchy problem
\begin{equation}\begin{split}\label{one-third-n}
\frac{\partial u(t,x)}{\partial t} & =\sum_{j=1}^{n-1} \frac{t^{1-j/n}}{\Gamma(j/n)} L_x^j f(x)
+ L_x^n u(t,x); u(0,x)=  f(x)\quad 
\end{split}\end{equation}
and the fractional Cauchy problem (\ref{frac-derivative-0}) with $\beta=1/n$
have the same unique solution given by
$
u(t,x) =\int_0^\infty p((t/s)^\beta,x) g_\beta(s)\,ds
$
with $\beta=1/n$.  
Hence the process $Z_t=X(E_t)$ is also the stochastic solution to this higher order Cauchy problem.


\section{\bf  Other subordinators}

An $\alpha$-time process is a Markov process subordinated to the absolute value of an independent one-dimensional symmetric
$\alpha$-stable process: $Z_t=B(|S_t|)$, where $B_t$ is a Markov process and $S_t$ is an independent symmetric $\alpha$-stable process both started at $0$. Let $Z_t^x =x+ Z_t$ the process started at $x$.

This process is self similar with index $1/2\alpha$ when the outer
process $X$ is a Brownian motion. In this case, Nane \cite{nane4} defined
the local time of this process and obtained laws of the iterated
logarithm for the local time for large time.

\subsection{\bf PDE-connection:}

\begin{theorem}[Nane (2005) \cite{nane3}]

Let $T(s)f(x)=E[f(X^{x}(s))]$ be the semigroup of the
continuous Markov process $X^{x}(t) $ and let $ L_x$ be its
generator. Let $\alpha=1$. Let $f$ be a bounded measurable
function in the domain of $ L_x$, with $D_{ij}f$ bounded and
H\"{o}lder continuous for all $1\leq i,\ j \leq n$. Then
$u(t,x)=E[f(Z^{x}_t)]$
solves
\begin{eqnarray}
\frac{\partial ^{2}}{\partial t^{2}}u(t,x)\ &\ = \ & -\frac{2L_x f(x)}{\pi t}\ - \ L_x^{2}u(t,x); \ u(0,x)  =\   f(x).\nonumber
\end{eqnarray}
\end{theorem}

 For  $\alpha=l/m \neq 1$ rational: the PDE is more
complicated since kernels of symmetric $\alpha$-stable processes
satisfy a higher order PDE:
$$
(\frac{\partial^{2}}{\partial
s^{2}})^{l}+(-1)^{l+1}\frac{\partial^{2m}}{\partial
t^{2m}})p_{t}^{\alpha}(0,s)=0.
$$
We also have to assume that  we can take the operator out of the integral. This is
valid for $\alpha=1/m$, $m=2,3,\cdots$, by a Lemma
in Nane \cite{nane3}.

\begin{theorem}[Nane (2005)\cite{nane3}]
Let $\alpha \in (0,2)$ be a rational with $\alpha=l/m$, where $l$ and $m$
are relatively prime. Let $ T(s)f(x)=E[f(X^{x}(s))]$ be the
semigroup of the continuous Markov process $X^{x}(t) $ and let
$L_x$ be its generator. Let $f$ be a bounded measurable
function in the domain of $L_x$, with $D^{\gamma}f$  bounded and
H\"{o}lder continuous for all multi index $\gamma$ such that $|\gamma |=2l$. Then
$u(t,x)=E[f(Z^{x}_{t})]$
solves
\begin{eqnarray}
 (-1)^{l+1}\frac{\partial ^{2m}}{\partial t^{2m}}u(t,x)&=&
-2\sum_{i=1}^{l}\left(
\frac{\partial^{2l-2i}}{\partial s^{2l-2i}}p_{t}^{\alpha}(0,s)|_{s=0}\right)L_x^{2i-1}f(x)\ 
  -\ L_x^{2l}u(t,x); \nonumber\\
u(0,x)&   =  & f(x).\nonumber
\end{eqnarray}
\end{theorem}
For some other connections of PDE's and iterated processes, see papers by Nane \cite{nane3} and Allouba and Zheng \cite{allouba1}, Allouba \cite{allouba2}, Baeumer et al. \cite{nane-m-b} and references therein.

\section{\bf  Open Problems }

\noindent {\bf Question 1.} Looking at the governing PDE for
subordinators other than Brownian motion, are there any fractional
in time PDE which has the same solution as the higher order PDE?

\noindent {\bf Question 2.} Are there PDE connections of the
iterated processes in bounded domain as the PDE connection of
Brownian motion in bounded domains?

\section{\bf Acknowledgements:} Author thanks Professor Mark M. Meerschaert and Professor Yimin Xiao for their help and discussions on the results in this paper. I also would like to thank Professor Anatoly N. Kochubei for providing the references for the initial apperance of  equation (\ref{frac-derivative-0}) in the literature.

\end{document}